\documentclass[a4paper, 11pt, twoside]{article}
\usepackage{amsmath}
\usepackage{amscd}
\usepackage{amssymb}
\usepackage{amsthm}
\usepackage{graphics}
\usepackage{indentfirst}
\usepackage{latexsym}
\usepackage{amsfonts}
\usepackage{multicol}
\usepackage{stmaryrd}
\usepackage{array}
\usepackage{pxfonts}

\newtheoremstyle{mattthm}{}{}{\itshape}{}{\bfseries}{.}{ }{}

\theoremstyle{mattthm}
\newtheorem{lemma}{Lemma}[section]

\newtheorem{thm}[lemma]{Theorem}
\newtheorem{cory}[lemma]{Corollary}

\newtheoremstyle{mattdef}{}{}{}{}{\bfseries}{.}{ }{}

\theoremstyle{mattdef}

\newtheorem*{acks}{Acknowledgements}

\begin{document}

\newenvironment{pf}{\noindent\textbf{Proof.}}{\hfill \qedsymbol\newline}
\newenvironment{pfof}[1]{\vspace{\topsep}\noindent\textbf{Proof of {#1}.}}{\hfill \qedsymbol\newline}
\newenvironment{pfenum}{\noindent\textbf{Proof.}\indent\begin{enumerate}\vspace{-\topsep}}{\end{enumerate}\vspace{-\topsep}\hfill \qedsymbol\newline}
\newenvironment{pfnb}{\noindent\textbf{Proof.}}{\newline}

\newlength\raiser
\newcommand\rh[1]{\raisebox{-3pt}{#1}}
\newcommand\mo{-\negthinspace1}
\newcommand\hdom{\linethickness{0.5mm}\put(0,0){\line(1,0){24}}\put(0,12){\line(1,0){24}}\put(0,0){\line(0,1){12}}\put(24,0){\line(0,1){12}}}
\newcommand\vdom{\linethickness{0.5mm}\put(0,0){\line(0,1){24}}\put(12,0){\line(0,1){24}}\put(0,0){\line(1,0){12}}\put(0,24){\line(1,0){12}}}
\newcommand\lra\longrightarrow
\newcommand\rez{\operatorname{res}}
\newcommand\bsm{\begin{smallmatrix}}
\newcommand\esm{\end{smallmatrix}}
\newcommand{\rt}[1]{\rotatebox{90}{$#1$}}
\newcommand\la\lambda
\newcommand{\ol}{\overline}
\newcommand{\ul}{\underline}
\newcommand\reg{^{\operatorname{reg}}}
\newcommand{\lan}{\langle}
\newcommand{\ran}{\rangle}
\newcommand\fkh{\mathfrak{h}}
\newcommand\fka{\mathfrak{a}}
\newcommand\partn{\mathcal{P}}
\newcommand{\py}[3]{\,_{#1}{#2}_{#3}}
\newcommand{\pyy}[5]{\,_{#1}{#2}_{#3}{#4}_{#5}}
\newcommand{\thmlc}[3]{\textup{\textbf{(\!\! #1 \cite[#3]{#2})}}}
\newcommand{\sss}{\mathfrak{S}_}
\newcommand{\dom}{\trianglerighteqslant}
\newcommand{\doms}{\vartriangleright}
\newcommand{\ndom}{\ntrianglerighteqslant}
\newcommand{\ndoms}{\not\vartriangleright}
\newcommand{\domby}{\trianglelefteqslant}
\newcommand{\domsby}{\vartriangleleft}
\newcommand{\ndomby}{\ntrianglelefteqslant}
\newcommand{\ndomsby}{\not\vartriangleleft}
\newcommand{\subs}[1]{\subsection{#1}}
\newcommand{\nin}{\notin}
\newcommand{\nchar}{\operatorname{char}}
\newcommand{\thmcite}[2]{\textup{\textbf{\cite[#2]{#1}}}\ }
\newcommand\zez{\mathbb{Z}/e\mathbb{Z}}
\newcommand\zepz{\mathbb{Z}/(e+1)\mathbb{Z}}
\newcommand{\bbf}{\mathbb{F}}
\newcommand{\bbc}{\mathbb{C}}
\newcommand{\bbn}{\mathbb{N}}
\newcommand{\bbq}{\mathbb{Q}}
\newcommand{\bbz}{\mathbb{Z}}
\newcommand\zo{\bbn_0}
\newcommand{\gs}{\geqslant}
\newcommand{\ls}{\leqslant}
\renewcommand\leq\ls
\renewcommand\geq\gs
\newcommand\dw{^\triangle}
\newcommand\wod{^\triangledown}
\newcommand{\hhh}{\mathcal{H}_}
\newcommand{\hh}{\mathcal{H}}
\newcommand\hsn{\hhh{\bbf,\mo}(\sss n)}
\newcommand\hfq{\hhh{\bbf,q}(\sss n)}
\newcommand{\sect}[1]{\section{#1}}
\newcommand{\ff}{\mathfrak{f}}
\newcommand{\fff}{\mathfrak{F}}
\newcommand\cf{\mathcal{F}}
\newcommand\fkn{\mathfrak{n}}
\newcommand\sx{x}
\newcommand\bra[1]{|#1\ran}
\newcommand\arb[1]{\widehat{\bra{#1}}}
\newcommand\foc[1]{\mathcal{F}_{#1}}
\newcommand{\clam}{\begin{description}\item[\hspace{\leftmargin}Claim.]}
\newcommand{\prof}{\item[\hspace{\leftmargin}Proof.]}
\newcommand{\malc}{\end{description}}
\newcommand\ppmod[1]{\ (\operatorname{mod}\ #1)}
\newcommand\wed\wedge
\newcommand\wede\barwedge
\newcommand\uu[1]{\,\begin{array}{|@{\,}c@{\,}|}\hline #1\\\hline\end{array}\,}
\newcommand{\ux}[1]{\operatorname{ht}_{#1}}
\newcommand\erim{\operatorname{rim}}
\newcommand\mire{\operatorname{rim}'}
\newcommand\mmod{\ \operatorname{Mod}}
\newcommand\call{\mathcal{L}}
\newcommand\calt{\mathcal{T}}
\newcommand\calu{\mathcal{U}}
\newcommand\calv{\mathcal{V}}
\newcommand\calf{\mathcal{F}}
\newcommand\cgs\succcurlyeq
\newcommand\cls\preccurlyeq
\newcommand\inc{\operatorname{asc}}
\newcommand\qbinom[2]{\left[\begin{smallmatrix}#1\\#2\end{smallmatrix}\right]}
\newcommand\dqbinom[2]{\left[\begin{array}c#1\\#2\end{array}\right]}
\newcommand\spe[1]{S^{#1}_{\bbf,\mo}}
\newcommand\smp[1]{D^{#1}_{\bbf,\mo}}
\newcommand\per[1]{M^{#1}_{\bbf,\mo}}
\newcommand\speq[1]{S^{#1}_{\bbf,q}}
\newcommand\smpq[1]{D^{#1}_{\bbf,q}}
\newcommand\perq[1]{M^{#1}_{\bbf,q}}
\newcommand\bx{\put(0,0){\line(0,1){12}}\put(0,12){\line(1,0){12}}\put(12,0){\line(0,1){12}}\put(0,0){\line(1,0){12}}}


\renewcommand{\thepage}{}
\title{On the irreducible Specht modules for Iwahori--Hecke algebras of type A with $q=-\negthinspace1$}
\author{Matthew Fayers\\
\normalsize Queen Mary, University of London, Mile End Road, London E1 4NS, U.K.}\date{}
\maketitle
\begin{center}
\end{center}
\markboth{Matthew Fayers}{On the irreducible Specht modules for Iwahori--Hecke algebras of type A with $q=-\negthinspace1$}
\pagestyle{myheadings}

\begin{abstract}
Let $p$ be a prime and $\bbf$ a field of characteristic $p$, and let $\hhh n$ denote the Iwahori--Hecke algebra of the symmetric group $\sss n$ over $\bbf$ at $q=\mo$.  We prove that there are only finitely many partitions $\la$ such that both $\la$ and $\la'$ are $2$-singular and the Specht module $S^\la$ for $\hhh{|\la|}$ is irreducible.
\end{abstract}

\sect{Introduction}

Suppose $\bbf$ is a field and $q$ a non-zero element of $\bbf$.  The \emph{Iwahori--Hecke algebra} $\hhh n$ of type $A$ over $\bbf$ with parameter $q$ is a finite-dimensional algebra which arises in various contexts, and whose representation theory closely resembles the representation theory of the symmetric group in prime characteristic.  An important class of modules for $\hhh n$ is the class of \emph{Specht modules}, and an interesting problem is to determine exactly which Specht modules are irreducible.  This problem has been solved in all cases except when $q=\mo$ and the characteristic of $\bbf$ is not $2$.  Various partial results are known for this case, together with a conjectured solution \cite[Conjecture 2.2]{fl} for the case where $\bbf$ has infinite characteristic (we adopt the convention that the \emph{characteristic} of a field is the order of its prime subfield).  In this note, we concentrate on the case of finite characteristic $p$.  Since the reducibility or not of Specht modules labelled by $2$-regular partitions and their conjugates is known, we can concentrate on partitions $\la$ such that neither $\la$ nor $\la'$ is $2$-regular.  Our main result is that for each $p$ there are only finitely many such partitions which label irreducible Specht modules.  Our approach is to use a decomposition map to relate the problem to the (known) classification of Specht modules at a $2p$th root of unity in infinite characteristic, and then to apply a recent result of the author and Lyle which proves the reducibility of a large class of Specht modules.  We complete the proof by employing some simple combinatorics of partitions.

The next section contains the necessary background material, and a statement and proof of the main result.  In Section \ref{further}, we try to examine more precisely the reducibility of Specht modules in prime characteristic, and give some computational results for small primes. 



\begin{acks}
The author would like to thank Sin\' ead Lyle for helpful discussions on this subject.
\end{acks}

\sect{The main result}\label{mainthmsec}

Suppose $\bbf$ is a field, and $q$ is a non-zero element of $\bbf$; we define $e=e(q)$ to be the multiplicative order of $q$ in $\bbf$ if $q\neq1$, or $e=\nchar(\bbf)$ if $q=1$.  For any $n\gs0$, we define $\hhh n=\hhh{\bbf,q}(\sss n)$ to be the Iwahori--Hecke algebra of the symmetric group $\sss n$ with parameter $q$.  The essential reference for the representation theory of $\hhh n$ is Mathas's book \cite{mathbook}.

Many of the important representations of $\hhh n$ are labelled by \emph{partitions} of $n$.  Recall that a partition of $n$ is a weakly decreasing sequence $\la=(\la_1,\la_2,\dots)$ of non-negative integers which sum to $n$.  When writing partitions, we usually group equal parts with a superscript and omit trailing zeroes.

For every partition $\la$ of $n$, there is an $\hhh n$-module $S^\la$ called the \emph{Specht module}.  (Note that we refer to the Specht module defined by Dipper and James \cite{dija}, rather than that used by Mathas.)  In the case $e=\infty$, the Specht modules are irreducible and afford all the irreducible representations of $\hhh n$.  In the case where $e$ is finite, the Specht modules may be reducible, but the irreducible $\hhh n$-modules can be obtained from them.  Let us say that a partition $\la$ is \emph{$e$-regular} if there is no $i$ such that $\la_i=\dots=\la_{i+e-1}>0$, and \emph{$e$-singular} otherwise.  When $\la$ is $e$-regular, the Specht module $S^\la$ has a unique irreducible quotient $D^\la$, and the modules $D^\la$ afford all the irreducible representations of $\hhh n$ as $\la$ ranges over the $e$-regular partitions of $n$.  The \emph{decomposition matrix} of $\hhh n$ has rows indexed by the partitions of $n$ and columns indexed by the $e$-regular partitions of $n$, with the $(\la,\mu)$-entry being the composition multiplicity $[S^\la:D^\mu]$.  In the case where $p=\infty$, we denote this decomposition matrix $D^{(e)}_n$; it is known that (given $e$) this matrix does not depend on the choice of $q$.

The subject of this paper is the problem of classifying the irreducible Specht modules.  Let us say that a partition $\la$ is \emph{$(e,p)$-reducible} if the Specht module $S^\la$ is reducible when $\bbf$ has characteristic $p$, or \emph{$(e,p)$-irreducible} otherwise; this condition is known to depend only on $e$ and $p$, not on the particular choice of $\bbf$ and $q$.  In the case where $e>2$, the classification of $(e,p)$-reducible partitions has been completed in a series of papers \cite{jmjs,slred,mfreduc,mfirred,slcp,jlm}, but the case $e=2$ remains open.  This paper is a small contribution towards completing this case.

If $\la$ is a partition, let $\la'$ denote the conjugate partition, defined by
\[\la'_j = \left|\left\{i\ \left|\ \la_i\gs j\right.\right\}\right|.\]
It is known that $\la$ is $(e,p)$-reducible if and only if $\la'$ is; this is because $S^{\la'}$ is essentially the dual of $S^\la$ \cite[Exercise 3.14(iii)]{mathbook}.  Furthermore, if $\la$ is $e$-regular, it is known whether $\la$ is $(e,p)$-reducible \cite[Theorem 4.15]{jmjs}.  So in order to complete the classification of irreducible Specht modules, it suffices to consider the case $e=2$, and to consider only partitions $\la$ such that both $\la$ and $\la'$ are $2$-singular.  Let us say that $\la$ is \emph{doubly-singular} if this is the case.  Now we can state our main result.

\vspace{\topsep}

\noindent\hspace{-3pt}\fbox{\parbox{469pt}{\vspace{-\topsep}

\begin{thm}\label{main}
Suppose $p$ is a prime.  Then there are only finitely many doubly-singular $(2,p)$-irreducible partitions.
\end{thm}
\vspace{-\topsep}}}
\vspace{\topsep}

We remark that this theorem is certainly not true in the case $p=\infty$.  For example, any partition of the form $(a^b)$ is $(2,\infty)$-irreducible.  This was observed by Mathas, using \cite[Theorem 4.7]{jmq-1}.

To prove Theorem \ref{main}, we use the classification of $(2p,\infty)$-partitions and the theory of \emph{decomposition maps}.  An excellent introduction to decomposition maps can be found in Geck's article \cite{geck2}.  Using the set-up in Section 3 of \cite{geck2}, one can obtain the following (recall that $D^{(e)}_n$ denotes the decomposition matrix for an Iwahori--Hecke algebra at an $e$th root of unity in a field of infinite characteristic).

\begin{thm}\label{epiadj}
Suppose $p=\nchar(\bbf)$ is a prime and $e=2$, and let $D$ be the decomposition matrix of $\hhh n$.  Then for any non-negative integer $i$ there is a matrix $A$ with rows indexed by $(2p^i)$-regular partitions of $n$ and columns indexed by $2$-regular partitions of $n$, with the following properties:
\begin{itemize}
\item
the entries of $A$ are non-negative integers;
\item
there is at least one non-zero entry in each row of $A$;
\item
$D = D^{(2p^i)}_nA$.
\end{itemize}
\end{thm}

This result arises from a decomposition map between an Iwahori--Hecke algebra at a $(2p^i)$th root of unity in a field of infinite characteristic, and $\hhh n$.  The matrix $A$ is simply the decomposition matrix associated to this map.

As a consequence, we get the following.

\begin{cory}\label{22p}
Suppose $\la$ is a partition of $n$, and suppose that $\la$ is $(2p^i,\infty)$-reducible for some $i$.  Then $\la$ is $(2,p)$-reducible.
\end{cory}

\begin{pf}
Since $\la$ is $(2p^i,\infty)$-reducible, the sum of the entries in the $\la$-row of $D^{(2p^i)}_n$ is at least $2$.  Now the properties of $A$ guarantee that the sum of the entries in the $\la$-row of $D$ is at least $2$, so that $\la$ is $(2,p)$-reducible.
\end{pf}

We shall use this result mainly in the case $i=1$, employing the known classification of $(2p,\infty)$-reducible partitions.  This is most easily stated in terms of hook lengths in the Young diagram.  Given a partition $\la$, recall that the \emph{Young diagram} $[\la]$ is the set
\[\left\{(i,j)\in\bbn^2\ \left|\ j\ls\la_i\right.\right\},\]
whose elements we call the \emph{nodes} of $\la$.  Given such a node $(i,j)$, define the \emph{hook length} $h_\la(i,j)$ to be the integer $1+\la_i-j+\la'_j-i$.  If the Young diagram is drawn with the English convention, this is the number of nodes of $\la$ directly below or directly to the right of $(i,j)$, including $(i,j)$ itself.

Now given $e>2$, say that $\la$ is an \emph{$e$-JM partition} if the following condition holds: for every $(i,j)\in[\la]$ for which $e$ divides $h_\la(i,j)$, we have either
\begin{itemize}
\item
$e$ divides $h_\la(i,k)$ for all $1\ls k\ls \la_i$, or
\item
$e$ divides $h_\la(k,j)$ for all $1\ls k\ls \la'_j$.
\end{itemize}

Then the following is a special case of the results in \cite{jmjs,mfirred,slcp}.

\begin{thm}\label{jm}
Suppose $e>2$.  Then a partition $\la$ is $(e,\infty)$-irreducible if and only if $\la$ is an $e$-JM partition.
\end{thm}

So in trying to classify $(2,p)$-irreducible partitions, we can restrict attention to $(2p)$-JM partitions.  We can also make another strong restriction, thanks to a recent result of the author and Lyle.  Say that a partition $\la$ is \emph{broken} if there exist $1<c<d$ such that $\la_{c-1}-\la_c>1$ and $\la_{d-1}=\la_d>0$, and \emph{unbroken} otherwise.  Note that if $\la$ is broken, then so is $\la'$.

\begin{thm}\thmcite{fl}{Theorem 2.1}\label{flbroken}
Suppose $\la$ is a broken partition.  Then $\la$ is $(2,p)$-reducible.
\end{thm}

Applying Corollary \ref{22p}, we find that any doubly-singular partition $\la$ which is $(2,p)$-irreducible must be an unbroken $(2p)$-JM partition.  So in order to prove Theorem \ref{main}, it suffices to show that there are only finitely many such partitions, for any $p$.  This follows from a few simple combinatorial results.

Given a partition $\la$ and an integer $e>1$, say that $\la$ is an \emph{$e$-core} if there is no $(i,j)\in[\la]$ for which $e$ divides $h_\la(i,j)$.  Obviously, if $\la$ is an $e$-core, then $\la$ is an $e$-JM partition.

\begin{lemma}\label{jmcore}
Suppose $\la$ is an unbroken doubly-singular partition, and is a $(2p)$-JM partition.  Then $\la$ is a $2p$-core.
\end{lemma}

\begin{pf}
If not, then we have $2p\mid h_\la(i,j)$ for some $(i,j)\in[\la]$.  Since $\la$ is a $(2p)$-JM partition, we have either $2p\mid h_\la(i,k)$ for all $1\ls k\ls \la_i$, or $2p\mid h_\la(k,j)$ for all $1\ls k\ls\la'_j$.  By replacing $\la$ with $\la'$ if necessary, we can assume the latter case.
\clam
$\la_1-\la_2>1$.
\prof
First suppose $\la_2<j$. Then $\la'_j=1$, so $\la_1=h_\la(1,j)+j-1$.  Since $h_\la(1,j)$ is divisible by $2p$, it is at least $2p$, and hence $\la_1\gs j+2p-1$.  So $\la_1-\la_2\gs 2p$.

On the other hand, suppose $\la_2\gs j$.  Then $h_\la(1,j)$ and $h_\la(2,j)$ are both divisible by $2p$, and $h_\la(1,j)-h_\la(2,j)=\la_1-\la_2+1$.  So $\la_1-\la_2$ is congruent to $-1$ modulo $2p$, and in particular is at least $2p-1$.
\malc
But if $\la$ is a $2$-singular partition with $\la_1-\la_2>1$, then $\la$ is broken; contradiction.
\end{pf}

\begin{lemma}\label{singcore}
Suppose $\la$ is a $2$-singular partition, and let $a$ be maximal such that $\la_{a-1}=\la_a>0$.  If $\la$ is a $2p$-core, then $\la_a\ls2p-2$.
\end{lemma}

\begin{pf}
Suppose for a contradiction that $\la_a\gs 2p-1$, and consider the hook lengths $h_\la(a-1,j)$ and $h_\la(a,j)$ for $1\ls j\ls \la_a$.  Note that for any $j$ we have $h_\la(a-1,j)=h_\la(a,j)+1$.  Furthermore, since $a$ is maximal such that $\la_{a-1}=\la_a$, we have $\la'_{j}-\la'_{j+1}=0$ or $1$ for any $1\ls j\ls \la_a-1$, and hence either $h_\la(a,j)=h_\la(a,j+1)+1$ or $h_\la(a,j)=h_\la(a-1,j+1)+1$.  This implies that for any $1\ls k\ls\la_a$, the set
\[\left\{\left.h_\la(a-1,j)\ \right|\ k\ls j\ls \la_a\right\}\cup\left\{\left.h_\la(a,j)\ \right|\ k\ls j\ls \la_a\right\}\]
equals the interval $\{1,2,\dots,l\}$ for some $l\gs\la_a-k+2$.  Taking $k=1$, we find that the hook lengths of $\la$ include $2p$, a contradiction.
\end{pf}

\begin{cory}\label{boundcore}
Suppose $\la$ is an unbroken doubly-singular partition, and is a $2p$-core.  Then $\la_1,\la'_1\ls 4p-6$.  Hence there are only finitely many unbroken doubly-singular $2p$-cores.
\end{cory}

\begin{pf}
Let $a$ be maximal such that $\la_{a-1}=\la_a>0$, and let $b$ be maximal such that $\la'_{b-1}=\la'_b>0$.  Applying Lemma \ref{singcore} to $\la$ and $\la'$, we have $\la_a,\la'_b\ls2p-2$.

The fact that $\la$ is unbroken implies that $\la_{i}-\la_{i+1}\ls1$ for $i=1,\dots,a-1$; in particular, since we have $\la_{\la'_b}-\la_{\la'_b+1}>1$, we must have $a\ls\la'_b$.  We also deduce that $\la_1\ls\la_a+a-2$.

So
\begin{align*}
\la_1&\ls\la_a+a-2\\
&\ls 2p-2+\la'_b-2\\
&\ls 2p-2+2p-2-2\\
&=4p-6;
\end{align*}
replacing $\la$ with $\la'$, we also get $\la'_1\ls4p-6$.  It is clear that bounding $\la_1$ and $\la'_1$ above leaves only finitely many partitions, so there are only finitely many unbroken doubly-singular $2p$-cores.
\end{pf}

This completes the proof of Theorem \ref{main}.

\sect{Some small values of $p$}\label{further}

We have shown that every doubly-singular $(2,p)$-irreducible partition must be an unbroken $2p$-core, and that there are only finitely many such partitions.  But not every such partition is $(2,p)$-irreducible: there are doubly-singular unbroken partitions which are $(2,\infty)$-reducible (for example, the partition $(4^2,1)$, which is a $2p$-core provided $p>3$), and these partitions are $(2,p)$-reducible by Corollary \ref{22p} with $i=0$.  But it seems reasonable to conjecture that every doubly-singular unbroken $2p$-core which is $(2,\infty)$-irreducible is also $(2,p)$-irreducible.  However, the author has very little evidence for this.  We end this paper by summarising the information we have for small values of $p$.

\subsection*{$p=2$}

In the case $e=p=2$, $\hhh n$ is actually the group algebra of the symmetric group, and the irreducible Specht modules in this case have been classified by James and Mathas \cite{jm2}.  And our results verify this classification: there is only one unbroken doubly-singular $4$-core, namely the partition $(2^2)$, and this is indeed the only doubly-singular $(2,2)$-irreducible partition.

\subsection*{$p=3$}

There are ten unbroken doubly-singular $6$-cores, namely the partitions
\[(2^2),\ (3^2),\ (2^3),\ (3^2,1),\ (3,2^2),\ (4^2),\ (3^2,2),\ (2^4),\ (3^3),\ (4,3^2,1).\]

It is not difficult to verify that these partitions are all $(2,3)$-irreducible; the tables in \cite[Appendix B]{mathbook} deal with all except $(4,3^2,1)$, for which one may use the fact that $(3^3)$ is $(2,3)$-irreducible together with the Branching Rule.  So we have completed the classification of $(2,3)$-irreducible partitions.

\subsection*{$p=5$}

There are $227$ unbroken doubly-singular $10$-cores, of which $115$ are $(2,\infty)$-irreducible. (Note that even though a general classification of $(2,\infty)$-partitions is still unknown, any single partition can be checked using the LLT algorithm \cite{llt} and Ariki's Theorem \cite{ari}.)  The author has been unable to determine whether these partitions are all $(2,5)$-irreducible.  The first case where the $(2,5)$-irreducibility is difficult to determine is the partition $(6^2,5,4)$.

\end{document}